\documentclass[final,1p,times]{elsarticle}

\usepackage{amssymb}
 \usepackage{amsthm}
\usepackage{amscd}
\usepackage{amsmath}
\usepackage{amsfonts}
\usepackage{amssymb}
\usepackage{graphicx}
\usepackage{color}
\newtheorem{theorem}{Theorem}

\newtheorem{lemma}[theorem]{Lemma}
\newtheorem{corollary}[theorem]{Corollary}

\usepackage{mathrsfs}
\usepackage{titletoc}

\newcommand{\ra}{\rightarrow}
\newcommand{\p}{\partial}
\newcommand{\f}{\frac}

\newcommand{\be}{\begin{equation}}
\renewcommand{\ra}{\rightarrow}
\newcommand{\ee}{\end{equation}}
\newcommand{\bea}{\begin{eqnarray}}
\newcommand{\eea}{\end{eqnarray}}
\newcommand{\bna}{\begin{eqnarray*}}
\newcommand{\ena}{\end{eqnarray*}}

\renewcommand{\le}{\left}
\newcommand{\ri}{\right}

\journal{***}

\begin{document}

\begin{frontmatter}

\title{A remark on energy estimates concerning extremals for Trudinger-Moser inequalities on a disc}

\author{Yunyan Yang}
 \ead{yunyanyang@ruc.edu.cn}
 \address{Department of Mathematics,
Renmin University of China, Beijing 100872, P. R. China}

\begin{abstract}
  In this short note, we generalized an energy estimate due to Malchiodi-Martinazzi (J. Eur. Math. Soc. 16 (2014) 893-908) and
  Mancini-Martinazzi (Calc. Var. (2017) 56:94). As an application, we used it to reprove existence of extremals for Trudinger-Moser inequalities
  of Adimurthi-Druet type on the unit disc. Such existence problems in general cases had been considered by Yang (Trans. Amer. Math. Soc. 359 (2007) 5761-5776;
  J. Differential Equations 258 (2015) 3161-3193) and Lu-Yang (Discrete Contin. Dyn. Syst. 25 (2009) 963-979) by using another method.
\end{abstract}

\begin{keyword}
Trudinger-Moser inequality\sep extremal function\sep energy estimate\sep blow-up analysis

\MSC[2010] 35A01\sep 35B09

\end{keyword}

\end{frontmatter}

\section{Introduction}
Let $\Omega\subset\mathbb{R}^2$ be a smooth bounded domain, $W_0^{1,2}(\Omega)$ be the usual
 Sobolev space, and $(u_\epsilon)_{\epsilon>0}$ be a sequence of
functions in $W_0^{1,2}(\Omega)$ satisfying the following equation
\be\label{equation}\le\{\begin{array}{lll}-\Delta u_\epsilon=\gamma_\epsilon u_\epsilon e^{u_\epsilon^2}
&{\rm in}& \Omega\\[1.5ex] u_\epsilon>0&{\rm in}& \Omega\\[1.5ex]u_\epsilon=0&{\rm on}& \p\Omega
\end{array}\ri.\ee
for some real number $\gamma_\epsilon>0$. It is of significance to study the blow-up behavior of $u_\epsilon$,
which was initiated by Adimurthi-Struwe \cite{Adi-Stru} and developed by Druet \cite{Druet}.
Particularly if $(u_\epsilon)_{\epsilon>0}$ is a sequence of solutions to (\ref{equation}) with $\|\nabla u_\epsilon\|_2^2\leq C$
and $\sup_{\Omega}u_\epsilon\ra\infty$ as $\epsilon\ra 0$, then it was proved by Druet \cite{Druet} that
there exists finite points $x_1,\cdots,x_\ell\in \overline{\Omega}$ such that
up to a subsequence $\gamma_\epsilon\ra \gamma_0\in [0,\lambda_1(\Omega)]$,
$u_\epsilon\ra u_0$ strongly in $C^1_{\rm loc}(\overline{\Omega}\setminus\{x_1,\cdots,x_\ell\})$ and weakly in $W_0^{1,2}(\Omega)$, where $\lambda_1(\Omega)$ is
the first eigenvalue of $-\Delta$ with respect to Dirichlet boundary condition,
\be\label{u0}-\Delta u_0=\gamma_0u_0e^{u_0^2}\quad{\rm in}\quad \Omega,\ee
and $\|\nabla u_\epsilon\|_2^2=4\pi N+\|\nabla u_0\|_2^2+o_\epsilon(1)$ for some integer $N\geq 1$.
This phenomenon is called the energy quantization. Later such a problem for critical growth equations in the sense of Trudinger-Moser embedding
received extensive interests, see for examples Martinazzi \cite{Mart}, Martinazzi-Struwe \cite{Mar-Stru},
Lamm-Robert-Struwe \cite{L-R-S}, Adimurthi-Yang \cite{Adi-Yang-Calc} and Yang \cite{Yang-Calc}.
Returning to (\ref{u0}), one may ask whether $u_0$ and $\gamma_0$ can actually be nonzero, and whether $N$ can actually be larger than $1$.
In the case $\Omega=B_1$, the unit disc in $\mathbb{R}^2$, Malchiodi-Martinazzi \cite{M-Martinazzi} gave a negative answer to this question,
namely
\\[1.5ex]
\noindent{\bf Theorem A}. {\it Let $(u_\epsilon)_{\epsilon>0}$ be a sequence of solutions to (\ref{equation}). Then either $(i)$
$\gamma_\epsilon\ra \gamma_0\in [0,\lambda_1(B_1)]$, $u_\epsilon\ra u_0$ in $C^1(\overline{B}_1)$, where
$u_0$ solves (\ref{u0}); or $(ii)$ $\gamma_\epsilon\ra 0$, $\|\nabla u_\epsilon\|_2^2\ra 4\pi$, $u_k\ra 0$ weakly in $W_0^{1,2}(B_1)$
and strongly in $C^1_{\rm loc}(\overline{B}_1\setminus\{0\})$, and $|\nabla u_\epsilon|^2dx\rightharpoonup 4\pi\delta_0$,
$\gamma_\epsilon u_\epsilon^2e^{u_\epsilon^2}dx\rightharpoonup 4\pi\delta_0$ weakly in the sense of measures.}\\[1.5ex]
 Note that in Theorem A, $\|\nabla u_\epsilon\|_2^2$ is not assumed to be bounded. It allows Theorem A to be applied to locate
 the critical values of the Trudinger-Moser functional
 $$E(u)=\int_{B_1}(e^{u^2}-1)dx$$
 constrained to the manifold
 $M_\Lambda=\{u\in W_0^{1,2}(B_1):\|\nabla u\|_2^2=\Lambda\}$. For details we refer the reader to (\cite{M-Martinazzi}, Theorem 1),
 which also gives the existence of blow-up solutions of (\ref{equation}) explicitly. Related works are referred to  
 Adimurthi-Prashanth \cite{A-Prashanth} and del Pino-Musso-Ruf \cite{d-M-Ruf1,d-M-Ruf2}.

Now let us summarize Malchiodi-Martinazzi's proof of Theorem A as follows: Let $c_\epsilon =u_\epsilon(0)=\max_{B_1}u_\epsilon$. If $c_\epsilon$ is bounded, then
 $(i)$ holds by elliptic estimates. Hereafter it is assumed that $c_\epsilon\ra\infty$ as $\epsilon\ra 0$. Set $r_\epsilon>0$ such that
$r_\epsilon^2\gamma_\epsilon c_\epsilon^2e^{c_\epsilon^2}=4$. Firstly they employed a function sequence on $B_{r_\epsilon^{-1}}$ as
\be\label{scale}\varphi_\epsilon(x)=c_\epsilon(u_\epsilon(r_\epsilon x)-c_\epsilon).\ee
It follows from (\ref{equation}) that
\be\label{funct}\le\{\begin{array}{lll}
-\Delta\varphi_\epsilon=4(1+c_\epsilon^{-2}\varphi_\epsilon)e^{2\varphi_\epsilon+c_\epsilon^{-2}\varphi_\epsilon^2}&{\rm in}& B_{r_\epsilon^{-1}}
\\[1.5ex]\varphi_\epsilon(0)=\sup_{B_{r_\epsilon^{-1}}}\varphi_\epsilon=0.
\end{array}\ri.\ee
Elliptic estimates applied on (\ref{funct}) lead to
$$\varphi_\epsilon\ra \varphi_0=-\log(1+|x|^2)\quad{\rm in}\quad C^1_{\rm loc}(\mathbb{R}^2).$$
Secondly they defined another function sequence $w_\epsilon=c_\epsilon^2(\varphi_\epsilon-\varphi_0)$ and proved that
\be\label{w-e}w_\epsilon\ra w_0\quad{\rm in}\quad C^1_{\rm loc}(\mathbb{R}^2),\ee where
\be\label{bub}w_0(r)=\varphi_0(r)+\f{2r^2}{1+r^2}-\f{1}{2}\varphi_0^2(r)+\f{1-r^2}{1+r^2}\int_{1}^{1+r^2}\f{\log t}{1-t}dt.\ee
 Thirdly they wrote
 $\varphi_\epsilon=\varphi_0+c_\epsilon^{-2}w_0+c_\epsilon^{-4}\eta_\epsilon$.
 Using the ODE theory and the contraction mapping theorem, they estimated the decay of $\eta_\epsilon$ and used it to derive
 $(ii)$ of Theorem A. Further analysis on asymptotic behavior of $w_0$ would lead to
 \be\label{energy}\|\nabla u_\epsilon\|_2^2=4\pi+O(c_\epsilon^{-4}).\ee

 To understand the sign of the error term $O(c_\epsilon^{-4})$ in (\ref{energy}), Mancini-Martinazzi \cite{Mancini-Martinazzi}
 refined the  analysis in \cite{M-Martinazzi}. Precisely they expanded
 $\varphi_\epsilon=\varphi_0+c_\epsilon^{-2}w_0+c_\epsilon^{-4}\eta_0+c_\epsilon^{-6}\phi_\epsilon$.
 Again using the contraction mapping theorem, they proved that $\phi_\epsilon$ is uniformly bounded up to large scales.
 Moreover, by analyzing the asymptotic behavior of the functions $w_0$ and $\eta_0$, they concluded the following:\\[1.5ex]
 \noindent{\bf Theorem B}. {\it Let $(u_\epsilon)_{\epsilon>0}$ be a sequence of solutions to (\ref{equation}). Assume $c_\epsilon=u_\epsilon(0)=
 \max_{B_1}u_\epsilon\ra\infty$ as $\epsilon\ra 0$. Then there holds
 \be\label{ener-2}4\pi+\f{4\pi}{c_\epsilon^{4}}+o(c_\epsilon^{-4})\leq \|\nabla u_\epsilon\|_2^2\leq 4\pi+\f{6\pi}{c_\epsilon^{4}}
 +o(c_\epsilon^{-4}).\ee}
 $\quad$Surprisingly it was shown by Mancini-Martinazzi \cite{Mancini-Martinazzi} that the energy estimate (\ref{ener-2}) can be used to reprove the critical Trudinger-Moser inequality \cite{Moser-71} and the existence of
 its extremals \cite{C-C} on the disc. Namely there exists a function $u_0\in W_0^{1,2}(B_1)$ with $\|\nabla u_0\|_2=1$ such that
 \be\label{Moser}\int_{B_1}e^{4\pi u_0^2}dx=\sup_{u\in W_0^{1,2}(B_1),\,\|\nabla u\|_2\leq 1}\int_{B_1}e^{4\pi u^2}dx.\ee

 In this short note, we shall use the above method of energy estimates to reprove the existence of extremals for Trudinger-Moser inequalities
 of Adimurthi-Druet's type on the disc (see \cite{A-D,Tintarev,Yang-Tran-07,Lu-Yang-1,Yang-JDE-15} for general versions), analogous to (\ref{Moser}). To do this, instead of (\ref{equation}), we consider
 the equation
 \be\label{eqn-1}\le\{\begin{array}{lll}
-\Delta u_\epsilon-\alpha_\epsilon u_\epsilon=\gamma_\epsilon u_\epsilon e^{u_\epsilon^2}&{\rm in}& B_1\\
[1.5ex] u_\epsilon>0&{\rm in}& B_1\\[1.5ex]
u_\epsilon=0 &{\rm on}&\p B_1
\end{array}\ri.\ee
for two real numbers $\gamma_\epsilon$ and $\alpha_\epsilon$. Adapting the arguments of Malchiodi-Martinazzi \cite{M-Martinazzi}
and Mancini-Martinazzi \cite{Mancini-Martinazzi}, we shall prove the following:

\begin{theorem}\label{Thm1}
Let $(u_\epsilon)_{\epsilon>0}$ be a sequence of decreasing radially symmetric solutions to (\ref{eqn-1}),
where $\gamma_\epsilon>0$ and $\alpha_\epsilon\ra\alpha<\lambda_1(B_1)$.
 Let $c_\epsilon=u_\epsilon(0)=\sup_{B_1}u_\epsilon$.
 If $c_\epsilon\ra\infty$ and $\limsup_{\epsilon\ra 0}\|\nabla u_\epsilon\|_{2}<\infty$, then there holds
\be\label{energy-est}
4\pi+\f{4\pi}{c_\epsilon^{4}}+o(c_\epsilon^{-4})\leq\|\nabla u_\epsilon\|_{2}^2-\alpha_\epsilon\|u_\epsilon\|_{2}^2\leq 4\pi+\f{6\pi}{c_\epsilon^{4}}+o(c_\epsilon^{-4}).
\ee
\end{theorem}

Compared with Theorems A and B, there is an additional assumption $\limsup_{\epsilon\ra 0}\|\nabla u_\epsilon\|_{2}<\infty$ in our theorem.
This is sufficient to ensure the convergence of $\varphi_\epsilon$ and $w_\epsilon$ defined as in (\ref{scale}) and (\ref{w-e}).
As an application of Theorem \ref{Thm1}, we have
the following:
\begin{corollary}\label{cor-2}
Let $\alpha<\lambda_1(B_1)$ be fixed. There exists a function $u_\ast\in W_0^{1,2}(B_1)\cap C^1(\overline{B}_1)$ with $\|u_\ast\|_{1,\alpha}=1$ such that
$$\int_{B_1}e^{4\pi {u_\ast}^2}dx=\sup_{u\in W_0^{1,2}(B_1),\, \|u\|_{1,\alpha}\leq 1}\int_{B_1}e^{4\pi u^2}dx,$$
where $\|u\|_{1,\alpha}^2=\|\nabla u\|_2^2-\alpha\|u\|_2^2$.
\end{corollary}

This recovers (\cite{Yang-JDE-15}, Theorem 1) on the disc $B_1$. Another application of Theorem \ref{Thm1} reads

\begin{corollary}\label{cor-3}
There exists some $\alpha_0>0$ such that if $0\leq \alpha\leq\alpha_0$, then the supremum
\be\label{sup}\sup_{u\in W_0^{1,2}(B_1),\, \|\nabla u\|_{2}\leq 1}\int_{B_1}e^{4\pi u^2(1+\alpha\|u\|_{2}^2)}dx\ee
can be attained by some function $u^\ast\in W_0^{1,2}(B_1)\cap C^1(\overline{B}_1)$ with $\|\nabla u^\ast\|_2=1$.
\end{corollary}

Corollary \ref{cor-3} is a special case of (\cite{Lu-Yang-1}, Theorem 1.2). Recently using energy estimates, Mancini-Thizy \cite{Mancini-Thizy} proved
 nonexistence of extremals for the supremum (\ref{sup}) when $\alpha$ is close to $\lambda_1(\Omega)$, where $\Omega$ is a smooth bounded domain
in $\mathbb{R}^2$.\\

 In the remaining part of this note, we shall prove Theorem \ref{Thm1} and Corollaries \ref{cor-2} and \ref{cor-3}.
Throughout this note, we do not
distinguish sequence and subsequence, and various constants are often denoted by the same $C$.

\section{Proof of Theorem \ref{Thm1}}\label{Th1}
For the proof of Theorem \ref{Thm1}, we shall modify the arguments in \cite{M-Martinazzi,Mancini-Martinazzi}.
Let $(u_\epsilon)_{\epsilon>0}$ be a sequence of decreasing radially symmetric solutions to (\ref{eqn-1}) for some $\gamma_\epsilon>0$
and $\alpha_\epsilon\ra\alpha<\lambda_1(B_1)$,
$c_\epsilon=u_\epsilon(0)=\sup_{B_1}u_\epsilon\ra\infty$, $\|\nabla u_\epsilon\|_{2}^2\leq \Lambda$ and $r_\epsilon>0$ be such that
\be\label{sc-l}r_\epsilon^2c_\epsilon^2\gamma_\epsilon e^{c_\epsilon^2}=4.\ee

A key observation is the following:
\begin{lemma}\label{lemma-1}
There exists a constant $\beta_0>0$ such that $r_\epsilon^2e^{\beta_0 c_\epsilon^2}\ra 0$ as $\epsilon\ra 0$.
\end{lemma}
{\it Proof}.
In view of (\ref{eqn-1}), we have
$$\gamma_\epsilon=\f{\|u_\epsilon\|_{1,\alpha_\epsilon}^2}{\int_{B_1}u_\epsilon^2e^{u_\epsilon^2}dx},$$
where $\|u_\epsilon\|_{1,\alpha_\epsilon}^2=\|\nabla u_\epsilon\|_2^2-\alpha_\epsilon \|u_\epsilon\|_2^2$.
Since $\alpha_\epsilon\ra \alpha<\lambda_1(B_1)$, there exists a constant $C>0$ such that $C^{-1}\|\nabla u_\epsilon\|_2\leq
\|u_\epsilon\|_{1,\alpha_\epsilon}\leq C\|\nabla u_\epsilon\|_2$.
By the H\"older inequality, the Sobolev embedding theorem and (\ref{sc-l}), we have for any $\beta<1$
\be\label{r-e}r_\epsilon^2e^{\beta c_\epsilon^2}=\f{4}{c_\epsilon^2 e^{(1-\beta)c_\epsilon^2}\|u_\epsilon\|_{1,\alpha_\epsilon}^2}
\int_{B_1}u_\epsilon^2e^{u_\epsilon^2}dx\leq \f{4}{c_\epsilon^2\|u_\epsilon\|_{1,\alpha_\epsilon}^2}\int_{B_1}u_\epsilon^2
e^{\beta u_\epsilon^2}dx\leq \f{C}{c_\epsilon^2}\le(\int_{B_1}e^{2\beta u_\epsilon^2}dx\ri)^{1/2}\ee
for some constant $C$. Since $\|\nabla u_\epsilon\|_{2}^2\leq \Lambda$, if $\beta<\f{1}{2\Lambda}$, then
 Moser's inequality leads to
$$\int_{B_1}e^{2\beta u_\epsilon^2}dx\leq C.$$
This together with (\ref{r-e}) and the assumption $c_\epsilon\ra\infty$ gives the desired result. $\hfill\Box$\\

Let $\psi_\epsilon(x)=c_\epsilon^{-1}u_\epsilon(r_\epsilon x)$ and
$\varphi_\epsilon(x)=c_\epsilon(u_\epsilon(r_\epsilon x)-c_\epsilon)$ for $x\in B_{r_\epsilon^{-1}}$.
Then we have
$$\le\{\begin{array}{lll}
-\Delta\varphi_\epsilon=\alpha_\epsilon r_\epsilon^2c_\epsilon^2\psi_\epsilon+4(1+c_\epsilon^{-2}\varphi_\epsilon)
e^{2\varphi_\epsilon+c_\epsilon^{-2}\varphi_\epsilon^2}\\[1.5ex]
\varphi_\epsilon(0)=\sup\varphi_\epsilon=0.
\end{array}\ri.$$
By Lemma \ref{lemma-1}, $\alpha_\epsilon r_\epsilon^2c_\epsilon^2\psi_\epsilon$ converges to zero uniformly in
$x\in B_{r_\epsilon^{-1}}$.
As usual (see \cite{M-Martinazzi}, Lemma 3) we have that
$\varphi_\epsilon\ra \varphi_0$ in $C^1_{\rm loc}(\mathbb{R}^2)$, where $\varphi_0(x)=-\log(1+|x|^2)$, and that
$$\lim_{R\ra\infty}\lim_{\epsilon\ra 0}\int_{B_{Rr_\epsilon}}\gamma_\epsilon u_\epsilon^2e^{u_\epsilon^2}dx=\int_{\mathbb{R}^2}
4e^{2\varphi_0}dx=4\pi.$$
Let
$$w_\epsilon(x)=c_\epsilon^2(\varphi_\epsilon(x)-\varphi_0(x)),\quad x\in B_{r_\epsilon^{-1}}.$$
We calculate
\bna
-\Delta w_\epsilon=\alpha_\epsilon r_\epsilon^2c_\epsilon^4\psi_\epsilon+c_\epsilon^2\le(4(1+c_\epsilon^{-2}\varphi_\epsilon)
e^{(2+c_\epsilon^{-2}\varphi_\epsilon)\varphi_\epsilon}-4e^{2\varphi_0}\ri).
\ena
It follows from Lemma \ref{lemma-1} that $\alpha_\epsilon r_\epsilon^2c_\epsilon^4\psi_\epsilon\ra 0$ uniformly on $B_{r_\epsilon^{-1}}$.
Similar to (\ref{w-e}) (\cite{M-Martinazzi}, Lemma 4), there holds $w_\epsilon\ra w_0$ in $C^1_{\rm loc}(\mathbb{R}^2)$, where $w_0$ is given as in
(\ref{bub}).
Then repeating the argument of the proof of (\cite{Mancini-Martinazzi}, Propositions 12 and 14), we conclude
\bea\label{e-1}
&&\int_{B_{r_\epsilon e^{c_\epsilon}}}\gamma_\epsilon u_\epsilon^2e^{u_\epsilon^2}dx=4\pi+\f{4\pi}{c_\epsilon^2}+o(c_\epsilon^{-4}),\\\label{e-2}
&&\int_{B_1\setminus B_{r_\epsilon e^{c_\epsilon}}}\gamma_\epsilon u_\epsilon^2e^{u_\epsilon^2}dx\leq \f{2\pi}{c_\epsilon^4}+o(c_\epsilon^{-4}).
\eea
Note also that $r_\epsilon e^{c_\epsilon}\ra 0$ by the above Lemma \ref{lemma-1}.
Multiplying both sides of (\ref{eqn-1}) and integrating by parts, we obtain
$$\|u_\epsilon\|_{1,\alpha_\epsilon}^2=\int_{B_1}\gamma_\epsilon u_\epsilon^2e^{u_\epsilon^2}dx.$$
This together with (\ref{e-1}) and (\ref{e-2}) implies (\ref{energy-est}).
$\hfill\Box$

\section{Proofs of Corollaries \ref{cor-2} and \ref{cor-3}}\label{Cors}

\subsection{Proof of Corollary \ref{cor-2}}
Let $\alpha<\lambda_1(B_1)$ be fixed. According to (\cite{Yang-JDE-15}, the formula (18)),
for any $0<\epsilon<4\pi$, there exists a function $u_\epsilon\in W_0^{1,2}(B_1)\cap
C^1(\overline{B}_1)$ such that
$$\int_{B_1}e^{(4\pi-\epsilon)u_\epsilon^2}dx=\sup_{u\in W_0^{1,2}(B_1),\,\|u\|_{1,\alpha}\leq 1}
\int_{B_1}e^{(4\pi-\epsilon)u^2}dx.$$
Moreover, by using a symmetrization argument, we may assume $u_\epsilon$ is a decreasing radially symmetric function. The Euler-Lagrange
equation reads
\be\label{E-L}\le\{
  \begin{array}{lll}
  -\Delta u_\epsilon-\alpha u_\epsilon=\f{1}{\lambda_\epsilon}u_\epsilon e^{(4\pi-\epsilon)u_\epsilon^2}\,\,\,{\rm in}\,\,\,
  B_1,\\[1.5ex] u_\epsilon>0\,\,\,{\rm in}\,\,\,
  B_1,\\[1.5ex] u_\epsilon=0\,\,\,{\rm on}\,\,\,
  \p B_1,\\[1.5ex] \|u_\epsilon\|_{1,\alpha}=1,\\[1.5ex]
  \lambda_\epsilon=\int_{B_1} u_\epsilon^2 e^{(4\pi-\epsilon)u_\epsilon^2}dx.
  \end{array}
  \ri.\ee
  Let $\widetilde{u}_\epsilon=\sqrt{4\pi-\epsilon}\,u_\epsilon$. On one hand
  \be\label{sub}\|\widetilde{u}_\epsilon\|_{1,\alpha}^2=(4\pi-\epsilon)\|u_\epsilon\|_{1,\alpha}^2=4\pi-\epsilon.\ee
  On the other hand, if $c_\epsilon=u_\epsilon(0)\ra\infty$, then we have by using Theorem \ref{Thm1},
  $$\|\widetilde{u}_\epsilon\|_{1,\alpha}^2\geq 4\pi+4\pi c_\epsilon^{-4}+o(c_\epsilon^{-4}).$$
  This contradicts (\ref{sub}) when $\epsilon$ is sufficiently small and implies that $u_\epsilon$ is uniformly bounded in $B_1$. Then applying elliptic estimates to (\ref{E-L}),
  we have that $u_\epsilon$ converges to an extremal function $u_\ast$ in $C^1(\overline{B}_1)$, as desired. $\hfill\Box$

  \subsection{Proof of Corollary \ref{cor-3}}

  Using a direct method of variation and a symmetrization argument, for any $0<\epsilon<4\pi$ and
  $0\leq\alpha<\lambda_1(B_1)$, we can find a decreasing radially symmetric
  function $u_\epsilon\in W_0^{1,2}(B_1)\cap C^1(\overline{B}_1)$ satisfying
  $$\int_{B_1}e^{(4\pi-\epsilon)u_\epsilon^2(1+\alpha\|u_\epsilon\|_2^2)}dx=\sup_{u\in W_0^{1,2}(B_1),\,\|\nabla u\|_{2}\leq 1}
\int_{B_1}e^{(4\pi-\epsilon)u^2(1+\alpha\|u\|_2^2)}dx.$$
Moreover, $u_\epsilon$ is a solution of
\be\label{E-L-1}\le\{
  \begin{array}{lll}
  -\Delta u_\epsilon=\f{\beta_\epsilon}{\lambda_\epsilon}u_\epsilon e^{\alpha_\epsilon u_\epsilon^2}+\zeta_\epsilon u_\epsilon\,\,\,{\rm in}\,\,\,
  B_1\\[1.5ex] \|\nabla u_\epsilon\|_2=1,\quad u_\epsilon>0\,\,\,{\rm in}\,\,\,
  B_1\\[1.5ex] \alpha_\epsilon=(4\pi-\epsilon)(1+\alpha\|u_\epsilon\|_2^2)\\[1.5ex]
  \beta_\epsilon=(1+\alpha\|u_\epsilon\|_2^2)/(1+2\alpha\|u_\epsilon\|_2^2)\\[1.5ex]
  \zeta_\epsilon=\alpha/(1+2\alpha\|u_\epsilon\|_2^2)\\[1.5ex]
  \lambda_\epsilon=\int_{B_1} u_\epsilon^2 e^{(4\pi-\epsilon)u_\epsilon^2}dx.
  \end{array}
  \ri.\ee
 Suppose that $c_\epsilon=u_\epsilon(0)\ra\infty$. By (\cite{Lu-Yang-1}, Lemmas 3.1 and 3.6), $\|u_\epsilon\|_2\ra 0$ and $\|c_\epsilon u_\epsilon- G_\alpha\|_2\ra 0$ as $\epsilon\ra 0$, where $G_\alpha\in W_0^{1,q}(B_1)$ $(1<q<2)$ is a
 Green function satisfying $-\Delta G_\alpha=\delta_0+\alpha G_\alpha$. It then follows that
 \be\label{c-u}\|u_\epsilon\|_2^2=c_\epsilon^{-2}\|G_\alpha\|_2^2+o(c_\epsilon^{-2}).\ee
 Let $v_\epsilon=\sqrt{\alpha_\epsilon}\,u_\epsilon$. Since $\|\nabla u_\epsilon\|_2=1$ and $\|u_\epsilon\|_2\ra 0$, we have
 \be\label{s-b}\|\nabla v_\epsilon\|_2^2=\alpha_\epsilon=4\pi+o_\epsilon(1).\ee
 Clearly
 $$-\Delta v_\epsilon-\zeta_\epsilon v_\epsilon=\f{\beta_\epsilon}{\lambda_\epsilon}v_\epsilon e^{v_\epsilon^2}\quad{\rm in}\quad B_1.$$
 $\zeta_\epsilon\ra \alpha<\lambda_1(B_1)$ as $\epsilon\ra 0$. In view of (\ref{s-b}), there holds $\limsup_{\epsilon\ra 0}\|\nabla v_\epsilon\|_2=4\pi$.
 Applying Theorem \ref{Thm1}, we obtain
  $$\|\nabla v_\epsilon\|_2^2-\zeta_\epsilon\|v_\epsilon\|_2^2\geq 4\pi+\f{4\pi}{v_\epsilon^4(0)}+o(v_\epsilon^{-4}(0)). $$
 This leads to
 \be\label{3}\|\nabla u_\epsilon\|_2^2\geq \zeta_\epsilon\|u_\epsilon\|_2^2+\f{1}{\alpha_\epsilon}\le(4\pi+\f{4\pi}{\alpha_\epsilon c_\epsilon^4}+
 o(c_\epsilon^{-4})\ri).\ee
 One calculates
 \be\label{4}\zeta_\epsilon\|u_\epsilon\|_2^2=\alpha\|u_\epsilon\|_2^2-2\alpha^2\|u_\epsilon\|_2^4+o(\|u_\epsilon\|_2^{4}),\ee
 \be\label{5}\f{1}{\alpha_\epsilon}=\f{1}{4\pi-\epsilon}\le(1-\alpha\|u_\epsilon\|_2^2+{\alpha^2}\|u_\epsilon\|_2^4+o(\|u_\epsilon\|_2^4)\ri).\ee
 Inserting (\ref{4}) and (\ref{5}) into (\ref{3}), we have
 $$\|\nabla u_\epsilon\|_2^2\geq 1-\alpha^2\|u_\epsilon\|_2^4+\f{1}{4\pi c_\epsilon^4}+
 o(\|u_\epsilon\|_2^4)+o(c_\epsilon^{-4}).$$
 This together with (\ref{c-u}) gives
 \be\label{ener-geq}\|\nabla u_\epsilon\|_2^2\geq 1+\le(\f{1}{4\pi}-\alpha^2\|G_\alpha\|_2^4\ri)\f{1}{c_\epsilon^4}+o(c_\epsilon^{-4}).\ee
 Using elliptic estimates, one can find a sufficiently small $\alpha_0>0$ such that if $0\leq\alpha\leq\alpha_0$, then
 $$\alpha^2\|G_\alpha\|_2^4\leq \f{1}{8\pi}.$$
 Hence in the case $0\leq\alpha\leq \alpha_0$, (\ref{ener-geq}) leads to
 $$\|\nabla u_\epsilon\|_2^2\geq 1+\f{1}{8\pi c_\epsilon^4}+o(c_\epsilon^{-4}).$$
 This contradicts $\|\nabla u_\epsilon\|_2=1$ and concludes that $c_\epsilon$ must be bounded. Then applying elliptic estimates to
 (\ref{E-L-1}), we get the desired extremal. $\hfill\Box$

\bigskip

{\bf Acknowledgements}. This work is partly supported by the National Science Foundation of China, Grant No.
  11471014.


\begin{thebibliography}{00}

\bibitem{A-D} A. Adimurthi, O. Druet, Blow-up analysis in dimension 2 and a sharp form of Moser-Trudinger inequality, Comm. Partial Differential Equations 29 (2004) 295-322.
\bibitem{A-Prashanth} Adimurthi, S. Prashanth, Failure of Palais-Smale condition and blow-up analysis for the critical exponent problem in $\mathbb{R}^2$,
Proc. Indian Acad. Sci. Math. Sci. 107 (1997) 283-317.

\bibitem{Adi-Stru} A. Adimurthi, M. Struwe, Global compactness
properties of semilinear elliptic equation with critical exponential
growth, J. Functional Analysis 175 (2000) 125-167.

\bibitem{Adi-Yang-Calc} A. Adimurthi, Y. Yang, Multibubble analysis on $N$-Laplace equation in $\mathbb{R}^N$, Calc. Var. 40 (2011) 1-14.

\bibitem{C-C} L. Carleson, A. Chang, On the existence
of an extremal function for an inequality of J. Moser, Bull. Sci.
Math. 110 (1986) 113-127.


\bibitem{d-M-Ruf1} M. del Pino, M. Musso, B. Ruf, New solutions for Trudinger-Moser critical equations in $\mathbb{R}^2$, J. Funct. Anal. 258
(2010) 421-457.

\bibitem{d-M-Ruf2} M. del Pino, M. Musso, B. Ruf, Byond the Trudinger-Moser supremum, Calc. Var. 44 (2012) 543-567.

\bibitem{Druet} O. Druet, Multibumps analysis in dimension 2, quantification of blow-up levels, Duke Math. J. 132 (2006) 217-269.

\bibitem{L-R-S} T. Lamm, F. Robert, M. Struwe, The heat flow with a
critical exponential nonlinearity, J. Functional Analysis 257 (2009)
2951-2998.


\bibitem{Lu-Yang-1} G. Lu, Y. Yang, The sharp constant and extremal functions for Moser-Trudinger inequalities involving $L^{p}$ norms, Discrete and Continuous Dynamical Systems 25 (2009), 963-979.

\bibitem{M-Martinazzi} A. Malchiodi, L. Martinazzi, Critical points of the Moser-Trudinger functional on a disk, J. Eur. Math. Soc. 16 (2014)
893-908.

\bibitem{Mancini-Martinazzi} G. Mancini, L. Martinazzi, The Moser-Trudinger inequality and its extremals on a disk via energy estimates, Calc. Var. (2017)
56:94,  https://doi.org/10.1007/s00526-017-1184-y.

\bibitem{Mancini-Thizy} G. Mancini, P. Thizy, Non-existence of extremals for the Adimurthi-Druet inequality, arXiv: 1711.05022.



\bibitem{Mart} L. Martinazzi, A threshold phenomenon for embeddings
of $H_0^m$ into Orlicz spaces, Calc. Var. Partial Differential
Equations 36 (2009) 493-506.

\bibitem{Mar-Stru} L. Martinazzi, M. Struwe, Quantization for an
elliptic equation of order $2m$ with critical exponential
non-linearity, Math Z. 270 (2012) 453-486.

\bibitem{Moser-71} J. Moser, A sharp form of an inequality by N. Trudinger, Indiana Univ. Math. J. 20 (1971) 1077-1092.



\bibitem{Tintarev} C. Tintarev, Trudinger-Moser inequality with remainder terms, J. Funct. Anal. 266 (2014) 55-66.


\bibitem{Yang-Tran-07} Y. Yang, A sharp form of the Moser-Trudinger inequality on a compact Riemannian surface,
Trans. Amer. Math. Soc. 359 (2007), 5761-5776.

\bibitem{Yang-JDE-15} Y. Yang, Extremal functions for Trudinger-Moser inequalities of Adimurthi-Druet type in dimension two,
J. Differential Equations 258 (2015) 3161-3193.

\bibitem{Yang-Calc} Y. Yang, Quantization for an elliptic equation with critical exponential growth on compact Riemannian surface
without boundary, Calc. Var. 53 (2015) 901-941.

\end{thebibliography}
\end{document}